\documentclass[12pt]{article} 
\usepackage{latexsym,amsfonts,amssymb,amsmath} 
\newtheorem{thm}{Theorem}[section] 
\newtheorem{crl}[thm]{Corollary} 
 
\newtheorem{lm}[thm]{Lemma}

\newtheorem{dfn}[thm]{Definition} 
\newenvironment{proof}{\noindent{\bf Proof.}}{\noindent$\Box$\par\medskip}

\newcommand{\ad}{\mathrm {ad}}

\newcommand{\sudda}[1]{}

\newcommand{\gr}{\mathrm {gr}}  
 
\newcommand{\spa}{\mathrm {span}}  
\begin{document} 
 
\date{} 
\title{Solvable Infinite Filiform Lie Algebras}   
\author {Clas L\"ofwall}
\maketitle
\begin{abstract} 
An infinite filiform Lie algebra $L$ is residually nilpotent and its
graded associated with respect to the lower central series has smallest
possible dimension in each degree but is still infinite. This means
that $\gr(L)$ is of dimension two in degree one and of dimension one
in all higher degrees. We prove that if $L$ is solvable, then already
$[L,L]$ is abelian. The isomorphism classes in this case are given in
\cite{bra}, but the proof is incomplete. We make the necessary
additional computations and restate the result in \cite{bra} when the ground field
is the complex numbers.
 
\end{abstract}

\section{Introduction} 
Infinite filiform Lie algebras have been
studied among others by Fialowski \cite{fia}, Millionshchikov
\cite{mil} and Shalev-Zelmanov \cite{zel}. They may be seen as
projective limits of finite dimensional filiform Lie algebras
introduced by Vergne \cite{ver} as nilpotent Lie algebras with maximal
degree of nilpotency among all nilpotent Lie algebras of a certain dimension. One result is that there is
only one infinite filiform naturally graded Lie algebra $L$, where
naturally graded means that $L$ is isomorphic to its graded associated
with respect to the filtration defined by the lower central series
($L^1=L,\ L^{i+1}=[L,L^i]$, $i\ge1$). This Lie algebra, denoted $M_0$, has a
basis $a,e_1,e_2,\ldots$ with $[a,e_i]=e_{i+1}$ and $[e_i,e_j]=0$ for
all $i$ and $j$. We have that $M_0$ is generated by $a,e_1$ which is a
basis for the component of degree
1, and $e_i$ is a basis for the component of degree $i$ for $i\ge2$. A
general infinite filiform Lie algebra  $L$ may be seen as a (filtered)
deformation of $M_0$ such that $\gr(L)\cong M_0$. Thus we have the
following definition:
\begin{dfn}An infinite filiform Lie algebra is a Lie algebra $L$ with 
 a basis $a,e_1,e_2,\ldots$ satisfying  
$$[a,e_j]=e_{j+1}\quad\text{and}\quad [e_j,e_{j+1}]=\sum_{s=1}^\infty \lambda_{js}e_{2j+1+s}\quad\text{for all}\quad
j\ge1
$$
for some $\lambda_{js}$, $j,s\ge1$ such that for each $j$,
$\lambda_{js}=0$ for $s>>1$
\end{dfn}
The missing products in the definition above are determined by the Jacobi identity, more
precisely, for $1\le k<m$ we have
$$
[e_k,e_m]=\sum_{s\ge1,0\le j-k\le m-j-1}(-1)^{j-k}\binom{m-j-1}{j-k}\lambda_{js}e_{m+k+s} 
$$
Moreover, the Jacobi identity impose conditions on $(\lambda_{js})$.
This is nicely described in \cite{mil}, where an explicit list of
quadratic polynomials are given with common zero-set equal to the set
of  $(\lambda_{js})$ that determines a deformation. The polynomials
are obtained by considering a deformation as above as a cochain $\psi$ of degree two
in the standard complex $C^*(M_0,M_0)$. The condition on $\psi$ making
the new structure defined by $\psi$ to a Lie algebra is in general that $\psi$ satisfies
the ``deformation equation", $d\psi+\frac{1}{2}[\psi,\psi]=0$, where
$[\cdot,\cdot]$ is the Nijenhuis-Richardson Lie superalgebra structure
on $C^*(M_0,M_0)$. It is proved in \cite{mil} that in this special
situation, the deformation equation is equivalent to
$d\psi=[\psi,\psi]=0$ and that  all possible
infinite filiform Lie algebras may be seen as the set of elements $\psi\in
\oplus_{i>0}H^{2,i}(M_0,M_0)$ such that $[\psi,\psi]=0$.

The solutions are known in two extreme cases. If $\lambda_{js}=0$ for
$j>1$ then there are no conditions on $(\lambda_{1s})$ to be a
solution. This gives rise to infinite filiform Lie algebras $L$
satisfying $[L,L]$ abelian, i.e., $L$ is solvable and next to abelian.
We will prove in the next section that these Lie algebras consist of
all solvable infinite filiform Lie algebras. The isomorphism classes,
which are parametrized by infinitely many parameters,
have been given by Bratzlavsky \cite{bra}, but the proof there is
incomplete. We perform the necessary extra computations in the last
section and restate Bratzlavsky's theorem. 

The other extreme case where the solutions are known is when $\lambda_{js}=0$ for
$s>1$. This gives only one new infinite filiform Lie algebra, namely
the subalgebra of the Witt algebra generated by
$x^{i+1}\partial/\partial x$ for $i\ge1$. 

The subalgebra of the Witt algebra generated by
$x^{2}d/dx$ and $x^{i+1}d/dx$ for $i\ge
k$ is also infinite filiform and defined by a deformation of $M_0$
satisfying $\lambda_{js}=0$ for $s\ne k-1$. It is conjectured in
\cite{mil} that this gives all solutions in the case when
$\lambda_{js}=0$ for $s\ne s_0$ (together with the solution when only
$\lambda_{1s_0}\ne0$).

We do not know if there are other infinite filiform Lie algebras than
those described above. In other words, does there exist a non-solvable infinite
filiform Lie algebra which is not a subalgebra of the Witt algebra?

\section{Criterion for solvability}
\begin{lm}\label{ab}Let $L$ be an infinite filiform Lie algebra. Then
$$
L \text{ is solvable } \Longleftrightarrow L^k \text{ is abelian for
some }k
$$
\end{lm}
\begin{proof} Let $L^{(1)}=L$ and $L^{(2k)}=[L^{(k)},L^{(k)}]$. Then
by induction and Jacobi identity $L^{(2k)}\subset L^{2k}$ and hence the implication to the
left follows. 

To prove the implication in the other direction, suppose $L$ is solvable and $L^2=[L,L]\ne0$. Then there is a
non-zero abelian ideal $I$ contained in $L^2$. Let $e_k+r\in I$, where $r\in L^{k+1}$. 
\begin{align*}
\text{\bf Claim:}\quad\quad\forall j\ge k\ \forall N \ \ e_j\in I+L^N
\end{align*}
Proof of Claim: By applying $\ad_a^{j-k}$ to $e_k+r\in I$ one gets
$e_j+r_1\in I$, where $r_1\in L^{j+1}$. Suppose $r_1=\lambda
e_{j+1}+r'$, where $r'\in L^{j+2}$. Since $\lambda
e_{j+1}+\lambda[a,r_1]\in I$, we get $e_j+r_2\in I$ for some $r_2\in
L^{j+2}$. Continuing in this manner, the claim follows. Now, let
$i,j\ge k$. Then for all $N\ge1$
$$
[e_i,e_j]\in [I,I]+L^N=L^N
$$
But $\cap_{N\ge1}L^N=0$. Hence $L^k$ is abelian.
\end{proof}
\begin{thm} \label{solv} Let $L$ be a solvable infinite filiform Lie algebra over a
field of characteristic zero. Then $[L,L]$ is abelian.
\end{thm}
\begin{proof} According to the lemma, we may 
suppose that $L^{n+1}$ is abelian and $n\ge2$. We will prove
that $L^n$ is abelian. This gives the theorem by induction. 

We know that $L$ has a basis $\{a,e_1,e_2,\ldots\}$ such that
$[a,e_i]=e_{i+1}$ and $[e_i,e_j]\in\spa\{e_k;\ k\ge
i+j+1\}$. We have $L^i=\spa\{e_j;\ j\ge i\}$. Hence, by assumption
$$
[e_i,e_j]=0\quad\text{for}\quad i,j\ge n+1
$$
Let
\begin{align*}
[e_{n-1},e_{n}]&=\lambda_0e_{2n}+\lambda_1e_{2n+1}+\ldots\\
[e_n,e_{n+1}]&=\mu_0e_{2n+2}+\mu_1e_{2n+3}+\ldots
\end{align*}
For $j>n+1$ we get
$$
[e_n,e_j]=[e_n,[a,e_{j-1}]]=-[e_{n+1},e_{j-1}]+[a,[e_n,e_{j-1}]]=[a,[e_n,e_{j-1}]]
$$
Hence, by induction
$$
[e_n,e_{j}]=\mu_0e_{n+j+1}+\mu_1e_{n+j+2}+\ldots
$$
for all $j\ge n+1$. Moreover,
\begin{align*}
[e_{n-1},e_{n+1}]&=[e_{n-1},[a,e_{n}]]=-[e_n,e_n]+[a,[e_{n-1},e_n]]\\
&=\lambda_0e_{2n+1}+\lambda_1e_{2n+2}+\ldots
\end{align*}
Suppose 
$$
[e_{n-1},e_{n+j}]=(\lambda_0-(j-1)\mu_0)e_{2n+j}+(\lambda_1-(j-1)\mu_1)e_{2n+j+1}+\ldots
$$
Then,
\begin{align*}
[e_{n-1},e_{n+j+1}]&=-[e_n,e_{n+j}]+[a,[e_{n-1},e_{n+j}]]\\
&=(\lambda_0-j\mu_0)e_{2n+j+1}+(\lambda_1-j\mu_1)e_{2n+j+2}+\ldots
\end{align*}
and hence, by induction, the last formula holds for all $j\ge0$. 
Since $L^{n+1}$ is abelian, we have
\begin{align*}
0&=[[e_{n-1},e_n],e_{n+1}]=\sum_{j\ge0}\mu_j[e_{n-1},e_{2n+2+j}]-\sum_{k\ge0}\lambda_k[e_n,e_{2n+1+k}]\\
&=\sum_{j\ge0}\mu_j\sum_{k\ge0}(\lambda_k-(n+j+1)\mu_k)e_{3n+j+2+k}-\sum_{j,k\ge0}\lambda_k\mu_je_{3n+k+2+j}\\
&=-\sum_{j,k\ge0}(n+j+1)\mu_j\mu_ke_{3n+2+j+k}
\end{align*}
Suppose $i\ge0$ and $\mu_j=0$ for $j\le i-1$. Then from above
$-(n+i+1)\mu_i^2=0$. Hence, $\mu_i=0$ and by induction it follows that
$\mu_i=0$ for all $i\ge0$. Hence $L^n$ is abelian.
\end{proof}
\begin{crl} Suppose a set $(\lambda_{js})$ of parameters determines a
deformation $L$ of $M_0$ as in the introduction and suppose $\lambda_{js}=0$ for all but finitely
many $j,s$. Then $\lambda_{js}=0$ for $j>1$.
\end{crl}
\begin{proof} The assumption gives that $L^k$ is abelian for some $k$. Hence,
by Lemma \ref{ab} and Theorem \ref{solv}, it follows that $[L,L]$ is abelian which
implies that $\lambda_{js}=0$ for $j>1$, since $[L,L]$ has a basis
$\{e_2,e_3,\ldots\}$ and $[e_j,e_{j+1}]=\sum_{s\ge1}\lambda_{js}e_{2j+1+s}$.  
\end{proof}
\section{The solvable case}
Suppose $(\lambda)=(\lambda_s)_{s\ge1}$ is an arbitrary sequence of
constants. We define a Lie algebra $L_{(\lambda)}$ as
follows. $L_{(\lambda)}$ has basis $\{a,e_1,e_2,\ldots\}$ and
multiplication
\begin{align*}[a,e_i]&=e_{i+1}\quad\text{for}\quad i\ge1\\
[e_1,e_i]&=\sum_{j\ge1}\lambda_je_{1+i+j}\quad\text{for}\quad i\ge2\\
[e_i,e_j]&=0\quad\text{for}\quad i,j\ge2
\end{align*}
It follows that $[a,[e_1,e_i]]=[e_1,e_{i+1}]$ for $i\ge2$. From this it is easy to
see that Jacobi identity holds, so $L_{(\lambda)}$ is indeed an
infinite filiform Lie
algebra such that $L_{(\lambda)}^2$ is abelian and any infinite filiform Lie
algebra $L$ such that $[L,L]$ is abelian is obtained in this way. Moreover, $L_{(\lambda)}$ is
obtained from $M_0$ by performing the deformation defined by
$(\lambda_{js})$, where $\lambda_{js}=0$ for $j>1$ and $\lambda_{1s}=\lambda_s$. 

In order to investigate when $L_{(\lambda)}$ and $L_{(\lambda')}$ are
isomorphic, we will study automorphisms $\phi$ of $L_{(\lambda)}$ and
determine the structure vector $(\lambda')$ in the new basis $\{\phi(a),\phi(e_1),\phi(e_2),\ldots\}$. Such a map
is determined
by
\begin{align*} \phi(a)&=c_0a+c_1e_1+c_2e_2+\ldots\\
\phi(e_1)&=d_0a+d_1e_1+d_2e_2+\ldots
\end{align*}
since then $\phi(e_{i+1})=\phi([a,e_i])=[\phi(a),\phi(e_i)]$ is
determined inductively. It is easily seen that $c_0\ne0$ and
$[\phi(e_1),\phi(e_2)]=c_0^{-1}d_0\phi(e_3)+\ldots$. Hence,
$\{\phi(a),\phi(e_1),\phi(e_2),\ldots\}$ is a basis of the same kind
as $\{a,e_1,e_2,\ldots\}$ (i.e., $[e_1,e_2]\in L^4$) only if
$d_0=0$. Then $\phi$ is an automorphism iff $c_0d_1\ne0$. Following \cite{gjk} we may
decompose any automorphism as a composition of three types of
automorphisms:
\begin{align*} \nu(c_0,c_1,d_1):&\quad\quad\begin{matrix}\phi(a)&=&c_0a&+&c_1e_1\\
\phi(e_1)&=&&&d_1e_1\end{matrix}\\
\\
\sigma(d,k),\ k\ge2:&\quad\quad\begin{matrix}\phi(a)&=&a&&\\
\phi(e_1)&=&e_1&+&de_k\end{matrix}\\
\\
\tau(c,k),\ k\ge2:&\quad\quad\begin{matrix}\phi(a)&=&a&+&ce_k\\
\phi(e_1)&=&e_1&&\end{matrix}
\end{align*}
The proof of the fact that the automorphisms of type $\sigma$ and
$\tau$ do not change $(\lambda)$ is missing in \cite{bra}. The
following lemma completes the proof in \cite{bra}. 
\begin{lm} Two Lie algebras $L_{(\lambda)}$ and $L_{(\lambda')}$ are
isomorphic iff $(\lambda')$ is obtained from $(\lambda)$ by performing
an automorphism of type $\nu$.
\end{lm} 
\begin{proof} Consider first the case $\sigma(d,k)$, $k\ge2$.  Then
\begin{align*}\phi(e_2)&=[a,e_1+de_k]=e_2+de_{k+1}\\
\phi(e_3)&=[a,e_2+de_{k+1}]=e_3+de_{k+2}
\end{align*}
and general, by induction
$$
\phi(e_i)=[a,e_{i-1}+de_{k+i-2}]=e_i+de_{k+i-1}
$$
Hence,
\begin{align*}[\phi(e_1),\phi(e_2)]=[e_1+de_k,e_2+de_{k+1}]=\sum_{i\ge1}\lambda_ie_{3+i}+d\sum_{i\ge1}\lambda_ie_{k+2+i}\\
=\sum_{i\ge1}\lambda_i(e_{3+i}+de_{k+2+i})=\sum_{i\ge1}\lambda_i\phi(e_{3+i})
\end{align*}
It follows that the automorphism of type $\sigma(d,k)$, $k\ge2$, does
not change the vector $(\lambda_{js})$.

Now, consider the case $\tau(c,k)$, $k\ge2$.
Then
\begin{align*}\phi(e_2)&=[a+ce_k,e_1]=e_2-c\sum_{i\ge1}\lambda_ie_{k+1+i}\\
\phi(e_3)&=[a+ce_k,e_2-c\sum_{i\ge1}\lambda_ie_{k+1+i}]=e_3-c\sum_{i\ge1}\lambda_ie_{k+2+i}
\end{align*}

and general, by induction
$$
\phi(e_j)=[a+ce_k,e_{j-1}-c\sum_{r\ge1}\lambda_re_{k+j-2+r}]=e_j-c\sum_{r\ge1}\lambda_re_{k+j-1+r}
$$
Hence,
\begin{align*}[\phi(e_1),\phi(e_2)]=[e_1,e_2-c\sum_{i\ge1}\lambda_ie_{k+1+i}]=\sum_{i\ge1}\lambda_ie_{3+i}-c\sum_{i\ge1}\lambda_i\sum_{r\ge1}\lambda_re_{k+2+i+r}\\
=\sum_{i\ge1}\lambda_i(e_{3+i}-c\sum_{r\ge1}\lambda_re_{k+2+i+r})=\sum_{i\ge1}\lambda_i\phi(e_{3+i})
\end{align*}
\end{proof} 
It is easy to see that an automorphism of type $\nu(c_0,0,d_1)$
transforms the sequence $(\lambda_1,\lambda_2,\ldots)$ into the sequence
$(c_0^{-2}d_1\lambda_1,c_0^{-3}d_1\lambda_2,\ldots)$. This gives the
possibility to choose the first two non-zero $\lambda_i:s$ to be $1$. 

Suppose that $\lambda_i=0$ for $i<t$ and $\lambda_t=1$ and consider
an automorphism $\phi$ of type $\nu(0,c,0)$. Then it is easy to see
that
\begin{align*}
[\phi(e_1),\phi(e_2)]=\phi(e_{3+t})+\lambda_{t+1}\phi(e_{3+t+1})+\lambda_{t+2}\phi(e_{3+t+2})+\ldots+\\
+(\lambda_{2t}-(1+t)c)\phi(e_{3+2t})+\
higher\ terms
\end{align*}
Hence if $c=\lambda_{2t}/(1+t)$ then $(\lambda)$ is transformed to
$(\lambda')$ where $\lambda'_i=0$ for $i<t$, $\lambda_t'=1$ and
$\lambda'_{2t}=0$. In all, any $(\lambda)$ may be transformed to
a sequence $(\mu)$ with the following properties
\begin{align*}
\mu_i&=\mu_j=1\quad\text{for some}\ 1\le i<j,\ j\ne2i\\
\mu_r&=0\quad\text{for}\ r<j,\ r\ne i\\
\mu_{2i}&=0
\end{align*}
Also, if $(\mu)$ and $(\mu')$ are two sequences of this kind, such that an
automorphism of type $\nu(c_0,c_1,d_1)$ transforms $(\mu)$ to $(\mu')$,
then this forces $c_0=1,c_1=0,d_1=1$ and hence $\mu=\mu'$.

Hence, we get the following version of
Bratzlavsky's theorem (in combination with Theorem \ref{solv}) in the
case when the ground
field is the complex numbers.
\begin{thm} Suppose $L$ is a solvable infinite filiform Lie
algebra over the complex numbers. Then there are unique integers $1\le r<s$, $s\ne2r$ and
complex numbers $\lambda_t$, $t>s$, $t\ne2r$ such that $L$ is
isomorphic to the infinite filiform Lie algebra given by the equations
\begin{align*}
[a,e_i]&=e_{i+1}\quad\text{for}\quad i\ge1\\
[e_1,e_i]&=e_{1+i+r}+e_{1+i+s}+\sum_{t>s,\
t\ne2r}\lambda_te_{1+i+t}\quad\text{for}\quad i\ge2\\
[e_i,e_{j}]&=0\quad \text{for}\quad i,j\ge2
\end{align*}  
\end{thm}

\end{document}